\documentclass[a4paper,11pt]{article}
\usepackage{amsfonts}
\usepackage{amsmath}
\usepackage{hyperref}

\input{tcilatex}
\begin{document}

\title{An Euler-type method for Volterra integro-differential equations}
\author{J. S. C. Prentice \\
Faculty of Information Technology\\
Eduvos\\
Bedfordview, South Africa\\
Email: justin.prentice@eduvos.com}
\maketitle

\begin{abstract}
We describe an algorithm, based on Euler's method, for solving Volterra
integro-differential equations. The algorithm approximates the relevant
integral by means of the composite Trapezium Rule, using the discrete nodes
of the independent variable as the required nodes for the integration
variable. We have developed an error control device, using Richardson
extrapolation, and we have achieved accuracy better than $10^{-12}$ for all
numerical examples considered.
\end{abstract}

\section{Introduction}

Many techniques exist for solving Volterra integro-differential equations
(IDEs), such as Adomian decomposition \cite{biazar}, Laplace decomposition 
\cite{bahuguna}, Galerkin methods \cite{maleknejad}, Haar functions \cite%
{malekjenad 2}, homotopy perturbation \cite{he 1}\ and more \cite{goghary}$-$%
\cite{hamoud}, including Runge-Kutta methods \cite{shimmary}\cite{brunner}.

In this \ paper, we focus our attention on Volterra IDEs of the form 
\begin{equation}
y^{\left( n\right) }\left( x\right) =f\left( x,y\right)
+\int\limits_{x_{0}}^{x}Kdt,\text{ \ \ \ }x>x_{0}  \label{problem}
\end{equation}%
with an appropriate set of initial conditions defined at $x_{0},$ and where
the kernel $K$ has the structure 
\begin{align}
(a)\text{ \ \ \ }K& =K\left( x,t\right)  \notag \\
(b)\text{ \ \ \ }K& =K\left( y\left( t\right) ,t\right)  \label{kernels} \\
(c)\text{ \ \ \ }K& =K\left( y^{\prime }\left( t\right) ,t\right)  \notag \\
(d)\text{ \ \ \ }K& =K_{1}\left( x\right) K_{2}\left( y\left( t\right)
,t\right) .  \notag
\end{align}%
The last of these is said to be \textit{separable}.

We will develop a straightforward one-step method, in the spirit of Euler,
which, combined with Richardson extrapolation, will be seen to yield very
accurate results.

Throughout this paper, we assume that all occurring functions are
real-valued and as smooth as our analysis requires.

\section{Algorithm}

Initially, we will describe our algorithm for the case of $n=1$ in (\ref%
{problem}). The more general case will be described later. We partition the
interval of interest, denoted $\left[ x_{0},x_{N}\right] ,$ by means of the
equispaced nodes%
\begin{equation}
x_{0}<x_{1}<x_{2}<\ldots <x_{N}.  \label{nodes}
\end{equation}%
The spacing between the nodes, known as the \textit{stepsize}, is denoted $h$%
. The stepsize must be constant in order for our error control device (based
on Richardson extrapolation) to be implemented successfully.

We assume that we have an initial value%
\begin{equation*}
y\left( x_{0}\right) =y_{0},
\end{equation*}%
and we compute the solution at $x_{1}$ via%
\begin{align*}
y_{1}& =y_{0}+hf\left( x_{0},y_{0}\right) +h\int\limits_{x_{0}}^{x_{0}}Kdt \\
& =y_{0}+hf\left( x_{0},y_{0}\right) .
\end{align*}%
This is an explicit Euler approximation to $y\left( x_{1}\right) .$

Then, we compute 
\begin{equation*}
y_{2}=y_{1}+hf\left( x_{1},y_{1}\right) +h\int\limits_{x_{0}}^{x_{1}}Kdt
\end{equation*}%
to obtain an approximation to $y\left( x_{2}\right) .$ Again, this step has
an explicit Eulerian character.

But how to find $\int\nolimits_{x_{0}}^{x_{1}}Kdt?$ To this end, we use the
information already determined, in the form%
\begin{equation*}
\int\limits_{x_{0}}^{x_{1}}Kdt\approx \frac{\left( x_{1}-x_{0}\right) }{2}%
\left( K_{0}+K_{1}\right) =\frac{h}{2}\left( K_{0}+K_{1}\right)
\end{equation*}%
where $K_{0}$ and $K_{1}$ denote the kernel $K$ evaluated at $x_{0}$ and $%
x_{1},$ respectively. This approximation is recognized as the Trapezium
Rule, wherein we have $t_{0}=x_{0}$ and $t_{1}=x_{1}$.

To find $y_{3},$ we compute%
\begin{align*}
y_{3}& =y_{2}+hf\left( x_{2},y_{2}\right) +h\int\limits_{x_{0}}^{x_{3}}Kdt \\
& =y_{2}+hf\left( x_{2},y_{2}\right) +\frac{h^{2}}{2}\left(
K_{0}+K_{1}\right) +\frac{h^{2}}{2}\left( K_{1}+K_{2}\right) \\
& =y_{2}+hf\left( x_{2},y_{2}\right) +\frac{h^{2}}{2}\left(
K_{0}+2K_{1}+K_{2}\right)
\end{align*}%
where the approximation to the integral is now seen to be the \textit{%
composite} Trapezium Rule, with $t_{0}=x_{0},t_{1}=x_{1}$ and $t_{2}=x_{2}.$

Continuing in this manner yields the general algorithm%
\begin{align*}
y_{i+1}& =y_{i}+hf\left( x_{i},y_{i}\right) +h\int\limits_{x_{0}}^{x_{i}}Kdt
\\
& =y_{i}+hf\left( x_{i},y_{i}\right) +\frac{h^{2}}{2}\left(
\sum\limits_{j=0}^{j=i}2K_{j}-\left( K_{0}+K_{i}\right) \right) .
\end{align*}

For the kernel $(c)$ in (\ref{kernels}), we simply express the derivative as%
\begin{equation*}
y_{j}^{\prime }=\frac{y_{j}-y_{j-1}}{h},
\end{equation*}%
and for kernel $(d)$, we have%
\begin{equation*}
y_{i+1}=y_{i}+hf\left( x_{i},y_{i}\right) +\frac{h^{2}K_{1}\left(
x_{i}\right) }{2}\left( \sum\limits_{j=0}^{j=i}2K_{2,j}-\left(
K_{2,0}+K_{2,i}\right) \right)
\end{equation*}%
i.e. we factor $K_{1}\left( x\right) $ out of the integral since it is not
dependent on $t$. For those kernels that are dependent on $y$ or $y^{\prime
} $, we have%
\begin{align*}
K\left( y\left( t_{j}\right) ,t_{j}\right) & =K\left( y_{j},x_{j}\right) \\
K\left( y^{\prime }\left( t_{j}\right) ,t_{j}\right) & =K\left(
y_{j}^{\prime },x_{j}\right) .
\end{align*}

When $n=2$ in (\ref{problem}), we have the system 
\begin{align*}
\left[ 
\begin{array}{c}
y^{\prime } \\ 
w^{\prime }%
\end{array}%
\right] & =\left[ 
\begin{array}{l}
w \\ 
f\left( x,y\right) +\int\limits_{x_{0}}^{x}Kdt%
\end{array}%
\right] \\
\Rightarrow \left[ 
\begin{array}{c}
y_{i+1} \\ 
w_{i+1}%
\end{array}%
\right] & =\left[ 
\begin{array}{l}
y_{i}+hw_{i} \\ 
w_{i}+hf\left( x_{i},y_{i}\right) +\frac{h^{2}}{2}\left(
\sum\limits_{j=0}^{j=i}2K_{j}-\left( K_{0}+K_{i}\right) \right)%
\end{array}%
\right]
\end{align*}%
and when $n=3,$ we have%
\begin{align*}
\left[ 
\begin{array}{c}
y^{\prime } \\ 
w^{\prime } \\ 
z^{\prime }%
\end{array}%
\right] & =\left[ 
\begin{array}{l}
w \\ 
z \\ 
f\left( x,y\right) +\int\limits_{x_{0}}^{x}Kdt%
\end{array}%
\right] \\
\Rightarrow \left[ 
\begin{array}{c}
y_{i+1} \\ 
w_{i+1} \\ 
z_{i+1}%
\end{array}%
\right] & =\left[ 
\begin{array}{l}
y_{i}+hw_{i} \\ 
w_{i}+hz_{i} \\ 
z_{i}+hf\left( x_{i},y_{i}\right) +\frac{h^{2}}{2}\left(
\sum\limits_{j=0}^{j=i}2K_{j}-\left( K_{0}+K_{i}\right) \right)%
\end{array}%
\right] .
\end{align*}%
Obviously, the initial values $y\left( x_{0}\right) $ and $w\left(
x_{0}\right) =y^{\prime }\left( x_{0}\right) $ must be specified for the
first system, and $y\left( x_{0}\right) ,w\left( x_{0}\right) =y^{\prime
}\left( x_{0}\right) $ and $z\left( x_{0}\right) =y^{\prime \prime }\left(
x_{0}\right) $ must be specified for the second system.

\section{Error control}

The Eulerian character of our algorithm, together with the use of the
Trapezium Rule, ensures that we cannot expect an error better than
first-order. However, this is quite acceptable, since we can deploy
Richardson extrapolation to achieve higher-order approximations from
first-order results. We have provided detail regarding Richardson
extrapolation elsewhere \cite{prentice 1}, and we simply state here the
process we use to construct solutions of order as high as five.

Let $y_{i}\left( h\right) $ denote the solution obtained at $x_{i}$ using a
stepsize $h$ (i.e. the nodes in (\ref{nodes})). Let $y_{i}\left( h/2\right) $
denote the solution obtained at $x_{i}$ using a stepsize $h/2$. Such a
computation uses the equispaced nodes%
\begin{equation*}
x_{0}<x_{1/2}<x_{1}<x_{3/2}<x_{2}<\ldots <x_{N-1/2}<x_{N}
\end{equation*}%
where each intermediate node $x_{i-1/2}$ is located midway between $x_{i-1}$
and $x_{i}$. We can similarly obtain the solutions $y_{i}\left( h/4\right)
,y_{i}\left( h/8\right) $ and $y_{i}\left( h/16\right) ,$ using appropriate
node distributions. Now, we form the linear combinations%
\begin{align*}
Y_{i}^{2}& =-y_{i}\left( h\right) +2y_{i}\left( h/2\right) \\
Y_{i}^{3}& =\frac{y_{i}\left( h\right) }{3}-2y_{i}\left( h/2\right) +\frac{%
8y_{i}\left( h/4\right) }{3} \\
Y_{i}^{4}& =-\frac{y_{i}\left( h\right) }{21}+\frac{2y_{i}\left( h/2\right) 
}{3}-\frac{8y_{i}\left( h/4\right) }{3}+\frac{64y_{i}\left( h/8\right) }{21}
\\
Y_{i}^{5}& =\frac{y_{i}\left( h\right) }{315}-\frac{2y_{i}\left( h/2\right) 
}{21}+\frac{8y_{i}\left( h/4\right) }{9}-\frac{64y_{i}\left( h/8\right) }{21}%
+\frac{1024y_{i}\left( h/16\right) }{315}
\end{align*}%
which yield 2nd-, 3rd-, 4th- and 5th-order solutions, respectively, at $%
x_{i} $. We will be interested in the 3rd-order solution in our numerical
examples. If we assume the 3rd- and 5th-order solutions have error terms of
the form%
\begin{align*}
& K_{3}h^{3}+\ldots \\
& K_{5}h^{5}+\ldots ,
\end{align*}%
respectively, then 
\begin{align*}
Y_{i}^{3}-Y_{i}^{5}& =K_{3}h^{3}+\ldots -\left( K_{5}h^{5}+\ldots \right) \\
& \approx K_{3}^{i}h^{3}
\end{align*}%
for suitably small $h$. Since $Y_{i}^{3}$ and $Y_{i}^{5}$ are known, we have 
\begin{equation*}
K_{3}^{i}=\frac{Y_{i}^{3}-Y_{i}^{5}}{h^{3}}
\end{equation*}%
as a good estimate for the error coefficient $K_{3}^{i}$. Consequently, a
suitable stepsize for a desired accuracy $\varepsilon $ is found from%
\begin{equation*}
h_{i}=\sigma \left( \frac{\varepsilon }{\left\vert K_{3}^{i}\right\vert }%
\right) ^{1/3}
\end{equation*}%
where the \textit{safety factor} $\sigma $ is $\sigma \sim 0.85.$ Naturally,
such a value for $h$ is computed at each $x_{i},$ and the smallest such
value is the one chosen. This chosen value is then used to rerun the
algorithm, with the resulting output satisfying the specified tolerance $%
\varepsilon .$ If we wish to control relative error, we compute%
\begin{equation*}
h_{i}=\sigma \left( \frac{\varepsilon \max \left\{ 1,\left\vert
y_{i}\right\vert \right\} }{\left\vert K_{3}\right\vert }\right) ^{1/3}
\end{equation*}%
at each $x_{i}$ and, as before, take the smallest such value and rerun the
algorithm.

\section{Examples}

\medskip We consider a variety of examples, indicated in the tables below.
For each example, we solve the IDE on the interval $\left[ 0,1\right] $ (see
the Appendix for commentary in this regard). The parameters $N_{1}$ and $%
N_{2}$ refer to the number of nodes $(N$ in (\ref{nodes})) needed to achieve
tolerances of $\varepsilon =10^{-6}$ and $\varepsilon =10^{-12},$
respectively, using the Richardson process described above. These examples
span the various possibilities in (\ref{problem}) and (\ref{kernels}). We\
have also included two examples of systems of IDEs (see Table 3). Initial
values used were determined from the given solutions, and so have not been
listed.

\begin{equation*}
\end{equation*}

\renewcommand{\arraystretch}{1.7}%
\begin{tabular}{|l|l|l|ll}
\multicolumn{4}{l}{Table 1: Examples $1-6,$ with values for $N_{1}$ and $%
N_{2}.$} & \multicolumn{1}{l}{} \\ \hline
\textbf{\#} & \textbf{IDE} & \textbf{Solution} & $N_{1}$ & 
\multicolumn{1}{|l|}{$N_{2}$} \\ \hline
\textbf{1} & $y^{\prime }=-1+\int_{0}^{x}y^{2}dt$ & $y\approx \frac{%
-x+x^{4}/28}{1+x^{3}/21}$ & $19$ & \multicolumn{1}{|l|}{$1892$} \\ \hline
\textbf{2} & $y^{\prime }=1+\int_{0}^{x}yy^{\prime }dt$ & $y=\sqrt{2}\tan
\left( x/\sqrt{2}\right) $ & $86$ & \multicolumn{1}{|l|}{$8513$} \\ \hline
\textbf{3} & $y^{\prime }=\cos x-\frac{x}{2}-\frac{\sin 2x}{4}%
+\int_{0}^{x}\left( y^{\prime }\right) ^{2}dt$ & $y=\sin x$ & $124$ & 
\multicolumn{1}{|l|}{$12322$} \\ \hline
\textbf{4} & $y^{\prime }=g\left( x\right) y-\int_{0}^{x}x^{2}t^{2}\left(
y^{\prime }\right) ^{3}dt$ & $y=\cos x$ & $38$ & \multicolumn{1}{|l|}{$3718$}
\\ \hline
\textbf{5} & $y^{\prime }=\left( -e^{x}-\frac{x^{2}e^{2x}}{3}\right)
y^{2}+\int_{0}^{x}\frac{t^{2}}{x}dt$ & $y=e^{-x}$ & $34$ & 
\multicolumn{1}{|l|}{$3344$} \\ \hline
\textbf{6} & $y^{\prime }=\frac{x^{2}+x+3}{3\left( x+1\right) }+\frac{%
2x^{3}-3x^{2}}{18}-\frac{y\left( x^{3}+1\right) }{3}+\int_{0}^{x}yt^{2}dt$ & 
$y=\ln \left( 1+x\right) $ & $17$ & \multicolumn{1}{|l|}{$1642$} \\ \hline
\end{tabular}%
\renewcommand{\arraystretch}{1}%
\begin{equation*}
\end{equation*}

The solution for \#1 is an approximation, as given in \cite{bahuguna}. In
\#4, we have%
\begin{equation*}
g\left( x\right) =\frac{\left( 
\begin{array}{c}
-27\sin x+27x^{4}\cos x-42x^{2}\cos x+2x^{2}\cos ^{3}x\ldots \text{ \ } \\ 
\text{ \ }-9x^{4}\cos ^{3}x-42x^{3}\sin x+6x^{3}\cos ^{2}x\sin x+40x^{2}%
\end{array}%
\right) }{27\cos x}.
\end{equation*}

\renewcommand{\arraystretch}{1.7}%
\begin{tabular}{|l|l|l|ll}
\multicolumn{4}{l}{Table 2: Examples $7-12,$ with values for $N_{1}$ and $%
N_{2}.$} &  \\ \hline
\textbf{\#} & \textbf{IDE} & \textbf{Solution} & $N_{1}$ & 
\multicolumn{1}{|l|}{$N_{2}$} \\ \hline
\textbf{7} & $y^{\prime }=3x^{2}-\frac{x^{4}}{3}+\int_{0}^{x}xt^{2}dt$ & $%
y=x^{3}$ & $26$ & \multicolumn{1}{|l|}{$2570$} \\ \hline
\textbf{8} & $y^{\prime }=y-\frac{x^{2}e^{x}}{2}+\int_{0}^{x}e^{x}tdt$ & $%
y=e^{x}$ & $45$ & \multicolumn{1}{|l|}{$4463$} \\ \hline
\textbf{9} & $y^{\prime }=\frac{2x^{3}+2x}{y+1}-\frac{x^{5}}{4}%
+\int_{0}^{x}xytdt$ & $y=x^{2}$ & $33$ & \multicolumn{1}{|l|}{$3238$} \\ 
\hline
\textbf{10} & $y^{\prime \prime }=x\cosh x-\int_{0}^{x}ytdt$ & $y=\sinh x$ & 
$30$ & \multicolumn{1}{|l|}{$2937$} \\ \hline
\textbf{11} & $y^{\prime \prime }=\left( \frac{\left( \ln \left( 1+x\right)
-1\right) \left( x+1\right) +1}{\left( x^{2}+1\right) \left(
4x^{2}+4x+1\right) }\right) y^{2}-\int_{0}^{x}\frac{\ln \left( t+1\right) }{%
x^{2}+1}dt$ & $y=2x+1$ & $13$ & \multicolumn{1}{|l|}{$1235$} \\ \hline
\textbf{12} & $y^{\prime \prime \prime }=e^{x}+e^{-x}-1+\int_{0}^{x}\frac{dt%
}{y}$ & $y=e^{x}$ & $33$ & \multicolumn{1}{|l|}{$3286$} \\ \hline
\end{tabular}%
\renewcommand{\arraystretch}{1}%
\begin{equation*}
\end{equation*}

\renewcommand{\arraystretch}{1.7}%
\begin{tabular}{|l|l|l|ll}
\multicolumn{4}{l}{Table 3: Examples $13-14,$ with values for $N_{1}$ and $%
N_{2}.$} &  \\ \hline
\textbf{\#} & \textbf{IDE} & \textbf{Solution} & $N_{1}$ & 
\multicolumn{1}{|l|}{$N_{2}$} \\ \hline
\textbf{13} & $%
\begin{array}{c}
y_{1}^{\prime }=2x-\frac{x^{5}}{5}-\frac{x^{10}}{10}+\int_{0}^{x}\left(
y_{1}^{2}+y_{2}^{3}\right) dt \\ 
y_{2}^{\prime }=3x^{2}+\int_{0}^{x}\left( y_{1}^{3}-y_{2}^{2}\right) dt\text{
\ \ \ \ \ \ \ \ \ \ \ \ \ }%
\end{array}%
$ & $%
\begin{array}{c}
y_{1}=x^{2} \\ 
y_{2}=x^{3}%
\end{array}%
$ & $83$ & \multicolumn{1}{|l|}{$8201$} \\ \hline
\textbf{14} & $%
\begin{array}{c}
y_{1}^{\prime }=1+x+x^{2}-y_{2}-\int_{0}^{x}\left( y_{1}+y_{2}\right) dt \\ 
y_{2}^{\prime }=-1-x+y_{1}-\int_{0}^{x}\left( y_{1}-y_{2}\right) dt\text{ \
\ \ }%
\end{array}%
$ & $%
\begin{array}{c}
y_{1}=x+e^{x} \\ 
y_{2}=x-e^{x}%
\end{array}%
$ & $30$ & \multicolumn{1}{|l|}{$2939$} \\ \hline
\end{tabular}%
\renewcommand{\arraystretch}{1}

\begin{equation*}
\end{equation*}

On our computational platform \cite{platform}, these calculations were
physically fast, requiring no more than five seconds, and usually much less,
for each case.

\section{Conclusion}

We have reported on an algorithm, based on Euler's method, for solving a
broad class of Volterra integro-differential equations. Our algorithm
approximates the relevant integral by means of the composite Trapezium Rule,
using the discrete nodes of the independent variable $x$ as the required
nodes for the integration variable $t$. We use Richardson extrapolation to
enhance the quality of the solution, achieving accuracy better than $%
10^{-12} $ for all the numerical examples considered. The algorithm has very
general character, is easy to implement and, on our computational platform,
is fast.

Nevertheless, further work is required. The algorithm is explicit, and we
have not considered stability issues in this work. It is possible that an
implicit form of the algorithm may be necessary to solve certain problems,
and the feasibility of such a version should be investigated. We believe
that for a \textit{nonseparable} kernel $K=K\left( x,y\left( t\right)
,y^{\prime }\left( t\right) ,t\right) $ a modification to the algorithm will
be necessary, and we will combine this task with that of creating an
implicit version. Lastly, we have not considered weakly singular problems
using our algorithm and this, too, should be a topic for further study.

\begin{equation*}
\end{equation*}

\begin{flushleft}
\textbf{{\Large {Appendix }}}

\ \ \ \ \ \ \ \ \ \ \ \ 
\end{flushleft}

Solving all the examples on $\left[ 0,1\right] $ is not restrictive.
Consider solving%
\begin{equation}
y^{\prime }\left( x\right) =f\left( x,y\left( x\right) \right)
+\int\limits_{x_{0}}^{x}K\left( x,y\left( t\right) ,t\right) dt
\label{original prob}
\end{equation}%
on $\left[ x_{0},x_{N}\right] .$ We can map $\left[ x_{0},x_{N}\right] $ to $%
\left[ 0,1\right] $ via%
\begin{equation*}
x=\left( x_{N}-x_{0}\right) \widetilde{x}+x_{0}\equiv m\widetilde{x}+x_{0},
\end{equation*}%
where $x\in \left[ x_{0},x_{N}\right] $ and $\widetilde{x}\in \left[ 0,1%
\right] .$ Hence,%
\begin{align*}
\frac{dy\left( x\right) }{dx}& =\frac{dy\left( m\widetilde{x}+x_{0}\right) }{%
d\widetilde{x}}\frac{d\widetilde{x}}{dx}=\left( \frac{1}{m}\right) \frac{%
dy\left( m\widetilde{x}+x_{0}\right) }{d\widetilde{x}} \\
f\left( x,y\left( x\right) \right) & =f\left( m\widetilde{x}+x_{0},y\left( m%
\widetilde{x}+x_{0}\right) \right) \\
K\left( x,y\left( t\right) ,t\right) & =K\left( m\widetilde{x}+x_{0},y\left(
m\widetilde{t}+x_{0}\right) ,m\widetilde{t}+x_{0}\right) \\
dt& =md\widetilde{t}
\end{align*}%
where the integration variable $t$ has been transformed in the same way as $%
x $.

We now have%
\begin{align}
\left( \frac{1}{m}\right) \frac{dy\left( m\widetilde{x}+x_{0}\right) }{d%
\widetilde{x}}=& f\left( m\widetilde{x}+x_{0},y\left( m\widetilde{x}%
+x_{0}\right) \right)  \notag \\
& +\int\limits_{0}^{\widetilde{x}}K\left( m\widetilde{x}+x_{0},y\left( m%
\widetilde{t}+x_{0}\right) ,m\widetilde{t}+x_{0}\right) md\widetilde{t} 
\notag \\
\Rightarrow \frac{d\widetilde{y}\left( \widetilde{x}\right) }{d\widetilde{x}}%
=& \widetilde{f}\left( \widetilde{x},\widetilde{y}\left( \widetilde{x}%
\right) \right) +\int\limits_{0}^{\widetilde{x}}\widetilde{K}\left( 
\widetilde{x},\widetilde{y}\left( \widetilde{t}\right) ,\widetilde{t}\right)
d\widetilde{t}  \label{dy tilde}
\end{align}%
where%
\begin{align*}
\widetilde{y}\left( \widetilde{x}\right) & \equiv y\left( m\widetilde{x}%
+x_{0}\right) \\
\widetilde{f}\left( \widetilde{x},\widetilde{y}\left( \widetilde{x}\right)
\right) & \equiv mf\left( m\widetilde{x}+x_{0},y\left( m\widetilde{x}%
+x_{0}\right) \right) \\
\widetilde{K}\left( \widetilde{x},\widetilde{y}\left( \widetilde{t}\right) ,%
\widetilde{t}\right) & \equiv m^{2}K\left( m\widetilde{x}+x_{0},y\left( m%
\widetilde{t}+x_{0}\right) ,m\widetilde{t}+x_{0}\right) .
\end{align*}

We solve (\ref{dy tilde}) for $\widetilde{y}\left( \widetilde{x}\right) $ on 
$\left[ 0,1\right] .$ The solution to the original problem (\ref{original
prob}) is then given by%
\begin{equation*}
y\left( x\right) =\widetilde{y}\left( \frac{x-x_{0}}{m}\right) .
\end{equation*}

For example, the IDE%
\begin{equation*}
y^{\prime }\left( x\right) =\left( \frac{-x^{5}+10x^{2}+32}{5x^{3}}\right)
y+\int\limits_{2}^{x}\frac{t^{2}y\left( t\right) }{x}dt,
\end{equation*}%
with $y\left( 2\right) =4,$ has the solution $y\left( x\right) =x^{2}.$ Say $%
\left[ x_{0},x_{N}\right] =\left[ 2,5\right] $ and we wish to transform the
problem to $\left[ 0,1\right] .$ We have%
\begin{equation*}
x=3\widetilde{x}+2\Rightarrow \widetilde{x}=\frac{x-2}{3}
\end{equation*}%
and $m=3$. Hence,%
\begin{align*}
\widetilde{f}\left( \widetilde{x},\widetilde{y}\left( \widetilde{x}\right)
\right) & =3\left( \frac{2}{3\widetilde{x}+2}-\frac{3\widetilde{x}+2}{4}%
+4\right) \widetilde{y} \\
\widetilde{K}\left( \widetilde{x},\widetilde{y}\left( \widetilde{t}\right) ,%
\widetilde{t}\right) & =9\left( \frac{\left( 3\widetilde{t}+2\right) ^{2}%
\widetilde{y}\left( \widetilde{t}\right) }{3\widetilde{x}+2}\right) .
\end{align*}%
It is easily confirmed that the solution to the transformed IDE%
\begin{equation*}
\frac{d\widetilde{y}\left( \widetilde{x}\right) }{d\widetilde{x}}=3\left( 
\frac{-\left( 3\widetilde{x}+2\right) ^{5}+10\left( 3\widetilde{x}+2\right)
^{2}+32}{5\left( 3\widetilde{x}+2\right) ^{3}}\right) \widetilde{y}%
+\int\limits_{0}^{\widetilde{x}}9\left( \frac{\left( 3\widetilde{t}+2\right)
^{2}\widetilde{y}\left( \widetilde{t}\right) }{3\widetilde{x}+2}\right) d%
\widetilde{t}
\end{equation*}%
is 
\begin{equation*}
\widetilde{y}\left( \widetilde{x}\right) =y\left( 3\widetilde{x}+2\right) =9%
\widetilde{x}^{2}+12\widetilde{x}+4.
\end{equation*}%
This gives%
\begin{align*}
y\left( x\right) & =\widetilde{y}\left( \frac{x-2}{3}\right) =9\left( \frac{%
x-2}{3}\right) ^{2}+12\left( \frac{x-2}{3}\right) +4 \\
& =x^{2}
\end{align*}%
as expected. Note that $\widetilde{y}\left( 0\right) =4=y\left( 2\right) .$

We see that an IDE defined on an arbitrary interval\ can be transformed to
the unit interval, and so we believe it is quite acceptable to solve all the
examples on $\left[ 0,1\right] .$

\end{document}